\documentclass[a4paper,10pt]{amsart}
\usepackage[english]{babel}
\usepackage{amssymb,amsmath,amscd,a4,amsfonts,amsthm,fullpage}
\date{\today}
\newtheoremstyle{style1}{10pt}{10pt}{\itshape}{}{\bfseries}{.}%
{7pt}{\thmnumber{\textbf{#2.}\;\;}\thmname{#1}\thmnote{#3}}
\newtheoremstyle{style2}{10pt}{10pt}{}{}{\bfseries}{.}%
{7pt}{\thmnumber{\textbf{#2.}\;\;}\thmname{#1}\thmnote{#3}}
\theoremstyle{style2}\newtheorem{dfn}{Definition}[section]
\theoremstyle{style2}
\theoremstyle{style2}
\theoremstyle{style2}\newtheorem{rem}[dfn]{Remark}
\theoremstyle{style2}\newtheorem{rems}[dfn]{Remarks}
\theoremstyle{style2}
\theoremstyle{style1}
\theoremstyle{style1}\newtheorem{prop}[dfn]{Proposition}
\theoremstyle{style1}\newtheorem{thm}[dfn]{Theorem}
\theoremstyle{style1}\newtheorem{cor}[dfn]{Corollary}

\begin{document}

\title{Some remarks on generalized roundness}

\author{Ghislain Jaudon}
\address{Universit\'{e} de Gen\`{e}ve,
Section de Math\'{e}matiques, 2-4 rue du Li\`{e}vre, Case postale
64, 1211 Gen\`{e}ve 4, Switzerland}
\email{ghislain.jaudon@math.unige.ch}

\begin{abstract}
\noindent By using the links between generalized roundness, negative type inequalities and equivariant Hilbert space compressions, we obtain that the generalized roundness of the usual Cayley graph of finitely generated free groups and free abelian groups of rank $\geq 2$ equals 1. This answers a question of J-F. Lafont and S. Prassidis.
\end{abstract}

\thanks{This work was supported by the
Swiss National Science Foundation Grant $\sharp$~PP002-68627.}

\maketitle

\setcounter{section}{0}
\section{Introduction}
\noindent Generalized roundness (see definition below) was introduced by P. Enflo in \cite{e1} and \cite{e2} in order to
study the uniform structure of metric spaces, and as an application
of this notion he gave a solution to Smirnov's problem
\cite{e2}. Rudiments of a general theory for generalized roundness
were developed in \cite{ltw}, where the link of this notion with
negative type inequalities is emphasized. More recently, generalized roundness was investigated in
the case of finitely generated groups \cite{lp}. Unfortunately,
generalized roundness is very difficult to estimate in general, and
for this reason there are only very few examples of metric spaces
for which the exact value is known. Here we use ideas developed in
\cite{ltw} together with estimates on generalized roundness
computed in \cite{lp} and results about equivariant Hilbert space
compression in \cite{gk} to deduce exact values of the generalized
roundness of finitely generated free (abelian and non-abelian)
groups endowed with their standard metrics.
\section{Preliminaries}
\noindent Let $(X,d)$ be a metric space, and let $G$ denote a group acting on $X$ by isometries.
\begin{dfn}
The \textit{generalized roundness of $(X,d)$} is the supremum of all positive numbers $p$ such that for every $n\geq 2$ and any collection of $2n$ points $\{a_1\ldots,a_n,b_1,\ldots,b_n\}$ in $X$, the following inequality holds:
$$
\sum_{1\leq i<j\leq n}(d(a_i,a_j)^p + d(b_i,b_j)^p) \leq \sum_{1\leq i,j\leq n} d(a_i,b_j)^p.
$$
We will denote the generalized roundness of the metric space $(X,d)$ by $\textrm{gr}(X,d)$, and simply $\textrm{gr}(X)$ when there is no ambiguity about the metric $d$.
\end{dfn}

\noindent Essentially, a metric space $(X,d)$ satisfies $\textrm{gr}(X,d)=p$ if $2n$-gons (for every $n\geq 2$) are thinner than the ones in $L^p$-spaces. This observation is justified by the following result (see \cite{ltw}):

\begin{prop}\label{gr de Lp}
Let $1 \leq p\leq 2$ and $(X,\mathcal{B},\mu)$ be a measured space. Then $\textrm{gr}(L^p(X,\mathcal{B},\mu))=p$.
\end{prop}

\begin{rem}\label{gr<=2}
The generalized roundness of any infinite and finitely generated group (endowed with the word metric) is always $\leq 2$ (see \cite{lp} Proposition 4.7).
\end{rem}

\begin{dfn}
A function $\psi:X\times X\rightarrow\mathbb{R}$
is said to be a {\it kernel of negative type} if $\psi(x,x)=0$ for all $x\in X$,
$\psi(x,y)=\psi(y,x)$ for all $x,y\in X$, and if for every integer
$n\geq 1$, for every $x_1,\ldots,x_n\in X$ and for every
$\lambda_1,\ldots,\lambda_n\in\mathbb{R}$ satisfying $\sum_{i=1}^n\lambda_i=0$, the following inequality
holds:
$$
\sum_{1\leq i,j\leq n}\lambda_i\lambda_j\psi(x_i,x_j)\leq 0.
$$
The kernel is said to be $G$-invariant if $\psi(gx,gy)=\psi(x,y)$ for all $x,y\in X$ and for all $g\in G$.
\end{dfn}

\noindent Kernels of negative type and generalized roundness are related by the following result (see \cite{ltw}):

\begin{thm}\label{rondeur vs type negatif}
$\textrm{gr}(X,d)\geq p$ if and only if $d^p$ is a kernel of negative type.
\end{thm}

\begin{dfn}
Let $\mathcal{H}$ be an Hilbert space. A map $f:X\rightarrow\mathcal{H}$ is said to be a \textit{uniform embedding} of $X$ into $\mathcal{H}$ if there exist non-decreasing functions $\rho_{\pm}(f):\mathbb{R}_+\rightarrow\mathbb{R}_+$ such that:
\begin{enumerate}
\item $\rho_-(f)(d(x,y))\leq \|f(x) - f(y)\|_{\mathcal{H}}\leq\rho_+(f)(d(x,y))$, for all $x,y\in X$;
\item $\lim_{r\rightarrow +\infty}\rho_{\pm}(f)(r)=+\infty$.
\end{enumerate}
Then the \textit{$G$-equivariant Hilbert space compression} of the
metric space $X$, denoted by $R_G(X)$, is defined as the supremum of
all $0<\beta\leq 1$ for which there exists a $G$-equivariant uniform
embedding $f$ into some Hilbert space which is equipped with an
action of $G$ by affine isometries, such that $\rho_+(f)$ is affine
and $\rho_-(f)(r)=r^{\beta}$ (for large enough $r$).
\end{dfn}

\noindent Concerning negative definite kernels, we will need a $G$-invariant analogue of the so-called GNS-construction (see for instance \cite{bhv}, 2.10):

\begin{prop}\label{GNS equivariant}
Let $\psi$ be a $G$-invariant kernel of negative type on $X$, then there exists a Hilbert space $\mathcal{H}$ equipped with an action of $G$ by affine isometries, and a $G$-equivariant map $f:X\rightarrow\mathcal{H}$, such that $\psi(x,y)=\|f(x) - f(y) \|_{\mathcal{H}}^2$ for all $x,y\in X$.
\end{prop}

\noindent Theorem \ref{rondeur vs type negatif} combined with Proposition \ref{GNS equivariant} immediately gives the following estimate:

\begin{prop}\label{inegalité}
For every group G of isometries of X, $R_G(X)\geq\frac{\textrm{gr}(X)}{2}$.
\end{prop}

\begin{rems}
On one hand, the previous inequality cannot be improved. Indeed, let
$X=G=\mathbb{Z}$ acting on itself by left translations and being endowed with its usual left invariant word metric.
Considering the inclusion of $\mathbb{Z}$ into $\mathbb{R}$, which is a $\mathbb{Z}$-equivariant isometry,
we have $R_{\mathbb{Z}}(\mathbb{Z})=1$ and moreover $\textrm{gr}(\mathbb{Z})\geq\textrm{gr}(\mathbb{R})$.
But Proposition \ref{gr de Lp} gives $\textrm{gr}(\mathbb{R})=2$. Therefore, by Remark \ref{gr<=2},
we obtain that $\textrm{gr}(\mathbb{Z})=2=2R_{\mathbb{Z}}(\mathbb{Z})$.\\
On the other hand, the inequality is unfortunately not an equality
in general. Consider for instance the case $X=G=\mathbb{Z}^2$. The
Hilbert space compression of $\mathbb{Z}^2$ equals 1 (see \cite{gk}
Example 2.7), and by amenability the equivariant Hilbert space
compression of $\mathbb{Z}^2$ equals the Hilbert space compression
(see \cite{ctv} Proposition 4.4). Hence
$R_{\mathbb{Z}^2}(\mathbb{Z}^2)=1$. But, by Corollary \ref{gr des
groupes libres} below, $\textrm{gr}(\mathbb{Z}^2)=1$.
\end{rems}

\section{Negative type inequalities in CAT(0) cube complexes}

\noindent Recall that a cube complex is a metric polyhedral complex
in which each cell is isometric to an Euclidean cube
$[-\frac{1}{2},\frac{1}{2}]^n$, and the gluing maps are isometries.
A finite dimensional cube complex always carries a complete geodesic
metric (see \cite{bh}). A cube complex is CAT(0) if it is simply
connected and if, in the link of every cube of the complex, there is
at most one edge between any two vertices and there is no triangle
not contained in a 2-simplex (see \cite{bh} and \cite{g}). Let $X$
denote a finite dimensional CAT(0) cube complex. The 0-skeleton
$X^{(0)}$ of $X$ can be endowed with the metric, denoted by $d_0$,
given by the length of the shortest edge path in the 1-skeleton of
$X$ between vertices. The proof of the next result is strongly
inspired by \cite{cn} Example 1.

\begin{thm}\label{gr pour cat0}
Let $X$ be a finite dimensional CAT(0) cube complex. Then $\textrm{gr}(X^{(0)}, d_0)\geq 1$.
\end{thm}

\begin{proof}
By Proposition \ref{gr de Lp}, it is sufficient to exhibit an
isometric embedding of $(X^{(0)}, d_0)$ into some $L^1$-space. Given
an edge in the complex, there is a unique isometrically embedded
codimension 1 coordinate hyperplane (again called hyperplane) which
cuts this edge transversely in its midpoint, and this hyperplane
separates the complex into two components, called half spaces (see
\cite{g}). We will denote by $H$ the set of all hyperplanes.
Moreover, by \cite{s}, shortest edge paths in the 1-skeleton cross
any hyperplane at most once. Hence the distance between two
vertices, $d_0(v,w)$, is the number of hyperplanes
separating $v$ and $w$ (hyperplanes such that the two vertices
are not in the same half space). We fix a vertex $v_0\in X^{(0)}$
and for every vertex $v\in X^{(0)}$ we set $H_v:=\{h\in H~\mid~
\textrm{h separates $v_0$ and $v$}\}$. Then we define
$$
f:X^{(0)}\rightarrow l^1(H), v\mapsto \sum_{h\in H_v}\delta_h
$$
where
$$
\delta_h:H\rightarrow\mathbb{R}, k\mapsto \begin{cases} 1~~~~ \textrm{if $k=h$}\\
                                                        0~~~~ \textrm{otherwise}
                                          \end{cases}
$$
It remains to show that $f$ is an isometry. Let $v,w$ be two vertices of $X$. Then
$$
\|f(v) - f(w)\|_{l^1(H)}=\sum_{l\in H}\left|\sum_{h\in H_v}\delta_h(l) - \sum_{h\in H_w}\delta_h(l)\right|.~~~~~~~~~~~~~~~~~~~~~~~~~~~~~~~~~~~~~~~~~~~~~~~~~~~~~~~~~~~~~~~~~~(\ast)
$$
For every $l\in H$, we have
$$
\left|\sum_{h\in H_v}\delta_h(l) - \sum_{h\in H_w}\delta_h(l)\right|=\begin{cases}  1 ~~~~ \textrm{if $l\in H_v\triangle H_w$}\\
                                                                       0 ~~~~ \textrm{if $l\in H_v\cap H_w$}
\end{cases}
$$
But a hyperplane $l\in H$ separates $v$ and $w$ if and only if $l\in H_v\triangle H_w$. Hence, the sum in the left member of $(\ast)$ is exactly the number of hyperplanes separating $v$ and $w$, i.e.,
$\|f(v) - f(w)\|_{l^1(H)}=d_0(v,w)$.

\end{proof}

\begin{cor}\label{gr des groupes libres} Let $n\geq 2$. We endow $\mathbb{Z}^n$ with the word metric associated to its canonical basis, and we endow the free group of rank $n$, $F_n$, with the word metric associated to any free generating system. We have:
\begin{enumerate}
\item $\textrm{gr}(\mathbb{Z}^n)=1$;

\item $\textrm{gr}(F_n)=1$.
\end{enumerate}

\end{cor}

\begin{proof}
$(i)$. By Corollary 4.14 of \cite{lp}, we have
$\textrm{gr}(\mathbb{Z}^n)\leq 1$. For the converse inequality, let
us consider the action of $\mathbb{Z}^n$ on $\mathbb{R}^n$ by left
translations. $\mathbb{R}^n$ can be viewed naturally as a CAT(0)
cube complex $X$ of which the 0-skeleton (endowed with the metric
$d_0$) is isometric to $\mathbb{Z}^n$.
Therefore, Theorem $\ref{gr pour cat0}$ gives the result.\\

\noindent $(ii)$. The Cayley graph of $F_n$ is a tree. In particular,
this is a 1-dimensional CAT(0) cube complex. Hence, by Theorem
\ref{gr pour cat0}, we obtain that $\textrm{gr}(F_n)\geq 1$. On the
other hand, it is known that $R_{F_n}(F_n)=\frac{1}{2}$ (see
\cite{gk}). Then by Proposition \ref{inegalité}, we deduce that
$\textrm{gr}(F_n)\leq 1$.
\end{proof}

\begin{rem}
Let $G$ be a group acting freely by isometries on the 0-skeleton
$(X^{(0)},d_0)$ of a CAT(0) cube complex $X$. If we fix a vertex
$v_0$, we define a metric $D_0$ on $G$ by setting
$D_0(g,h):=d_0(gv_0,hv_0)$, and Theorem \ref{gr pour cat0} gives
$\textrm{gr}(G,D_0)\geq 1$.
\end{rem}

\end{document}